\documentclass{svproc}
\usepackage{amsmath, verbatim, amsfonts, amssymb, color}
\usepackage{mathrsfs}
\usepackage{dsfont}
 
 \usepackage{graphicx}
 
\renewcommand\labelenumi{\textup{\alph{enumi})}}
\renewcommand\theenumi\labelenumi

 \makeatletter
\def\@makefnmark{\hbox{(\@textsuperscript{\normalfont\@thefnmark})}}
\makeatother

\makeatletter
\@namedef{subjclassname@2020}{%
  \textup{2020} Mathematics Subject Classification}
\makeatother

\def\M{{M}}%Riemannian manifold
\def\g{{g}}% tensor
\DeclareMathOperator{\ric}{\mathsf{Ric}}
\DeclareMathOperator{\Ric}{\mathrm{Ric}}

\DeclareMathOperator{\CD}{\mathsf{CD}}
\DeclareMathOperator{\RCD}{\mathsf{RCD}}
\DeclareMathOperator{\BE}{\mathsf{BE}}

\DeclareMathOperator{\Cpl}{\mathsf{Cpl}}

\def\Lip{\mathrm{Lip}}%Lip.stetige Funktionen

\def\N{{\mathbb N}}    %Zahlkoerper
\def\R{{\mathbb R}}%Zahlkoerper

\newcommand{\ee}{e}

\newcommand{\Pz}{\mathcal{P}}

\newcommand{\Ent}{\mbox{\rm Ent}}
\newcommand{\supp}{\mbox{\rm supp}}

\def\X{{\mathsf X}}
\def\d{{\mathsf d}}
\def\m{{\mathsf m}}
\def\M{{\mathsf M}}
\def\g{{\mathsf g}}
\def\Y{{\mathsf Y}}

\begin{document}
\mainmatter
\title{\bfseries Metric Measure Spaces and Synthetic Ricci Bounds\\ ---\\
Fundamental Concepts and Recent Developments}
\titlerunning{Metric Measure Spaces  and  Synthetic Ricci Bounds}

\author{Karl-Theodor Sturm}
\authorrunning{Sturm}
%\tocauthor{David Berger, Franziska K\"uhn and Ren\'e L. Schilling}
\institute{Hausdorff Center for Mathematics, University of Bonn, Germany \\
\email{sturm@uni-bonn.de}  }

\maketitle

\section{Synthetic Ricci Bounds for Metric Measure Spaces}

\subsection{Metric  Spaces %$(\X,\d)$
}

The class of \emph{metric spaces} $(\X,\d)$ is a far reaching 
generalizations of the class of \emph{Riemannian manifolds} $(\M,\g)$. It allows for rich geometric structures
including singularities, branching, change of dimension as well as fractional and infinite dimension.

Already in the middle of the last century, A.D.Alexandrov \cite{Alex48,Alex55}  has proposed his fundamental concepts of lower and upper bounds for 
generalized  sectional
curvature  for metric spaces. 
Especially these lower bounds are particularly well-behaved w.r.t.~the so-called Gromov-Hausdorff metric on the class
of compact metric spaces as observed 
by Gromov \cite{Gromov81,Gromov07}:
\begin{itemize}
\item[$\ast$] for each $K\in\R$, the class
$$\Big\{ (\X,\d) \mbox{ with sect. curv. }\ge K\Big\}$$ is closed
under GH-convergence;
\smallskip
\item[$\ast$] for each $K,L,N\in\R$, the class
$$\Big\{ (\X,\d) \mbox{ with sect. curv. }\ge K, \mbox{ dimension }\le N,
\mbox{ diameter }\le L \Big\}$$  is compact.
\end{itemize}
In the sequel, many properties of Riemannian manifolds and geometric estimates which only depend on one-sided curvature bounds  could be proven for such metric spaces $(\X,\d)$ with synthetic  (upper or lower) curvature bounds.  For  spaces   with synthetic lower bounds on the sectional curvature, also a far reaching analytic calculus was developed
with foundational contributions by
Burago--Gromov--Perelman \cite{BGP92},
Kuwae--Machigashira--Shioya \cite{KMS01},
Zhang--Zhu \cite{ZhangZhu12}.

\smallskip

However, for most properties and estimates in geometric analysis, spectral theory and stochastic analysis on manifolds, no quantitative assumptions on the sectional curvature are needed but --- as observed in the seminal works of Yau, Cheeger, Colding, Elworthy, Malliavin, Bismut, Perelman
and many others
--- merely a lower bound on the \emph{Ricci curvature}
$$\Ric\ge K\, \g.$$
Since the Ricci tensor is the trace of the sectional curvature, i.e.
$$\Ric_x(v_i,v_i):=\sum_{j\not=i}\operatorname{Sec}_x(v_i,v_j)\qquad \text{if $\{v_i\}_{i=1,\ldots,n}$ ONB of }T_xN,$$
assumptions on lower bounded Ricci curvature are less restrictive than assumptions on lower bounded sectional curvature.
Replacing (synthetic) sectional curvature bounds by Ricci bounds, the previously mentioned Gromov's compactness theorem turns into a precompactness theorem:
\begin{itemize}
\item[$\ast$] 
For any choice of $K,L,N\in\R$, the class of Riemannian manifolds $(\M,\g)$ with 
  Ricci curvature $\ge K$, dimension $\le N$ and diameter $\le L$ is relatively compact w.r.t.~mGH-convergence.
\end{itemize}
Properties of mGH-limits of Cauchy sequences in such classes (so-called \emph{Ricci limit spaces})  have been studied in great detail by Cheeger--Colding \cite{CC97,CC00a,CC00b}; see also \cite{ColdingNaber12,CheegerNaber13}.

As already pointed out by Gromov, the right setting to deal with the completions %closures or compactifications 
of these classes  is the class of \emph{metric measure spaces}.
However, what was missing for decades was a synthetic formulation of lower Ricci bounds, applicable not only to Riemannian manifolds (and their limits) but also to metric measure spaces.

\subsection{Metric Measure Spaces}
 
Here and in the sequel, a \emph{metric measure space} (briefly \emph{mm-space}) will always mean a triple $(\X,\d,\m)$ consisting of 
\begin{itemize}
\item a space $\X$,
\item a complete separable metric $\d$ on $\X$,
\item
a locally finite Borel  measure $\m$ on it. 
\end{itemize}
It is called normalized (or mm$_1$-space) iff in addition $\m(\X)=1$.

A primary goal since many years has been to find a formulation of
generalized  Ricci
curvature bounds $\mathrm{Ric}(\X,\d,\m)\ge K$ which is
\begin{itemize}
\item[$\ast$] equivalent to $\mathrm{Ric}_x(v,v)\ge K\,\|v\|^{2}$
if X is a Riemannian manifold,
\item[$\ast$] stable under convergence,
\item[$\ast$] intrinsic, synthetic (like curvature bounds in Alexandrov geometry),
\item[$\ast$] sufficient for many geometric, analytic, spectral theoretic conclusions.
\end{itemize}
In independent works, such a formulation   has been proposed by the author \cite{St06}, \cite{St06a}  and by Lott--Villani \cite{LV09}, based on the concept of optimal transport and relying on previous works by
 Brenier \cite{Brenier91}, Gangbo \cite{Gangbo94}, McCann  \cite{Mccann95,Mccann01}, Otto \cite{Otto01}, Otto--Villani \cite{OttoVill00}, Cordero-Erausquin--McCann--Schmuckenschl\"ager \cite{CMS01}, Sturm--vonRenesse \cite{StvRe05}.

The synthetic lower Ricci bound for a mm-space $(\X,\d,\m)$ will be  defined through the interplay of two quantities on $\X$:
\begin{itemize}
\item the \emph{Kantorovich--Wasserstein distance} 
\begin{equation}\label{W2}
W_2(\mu_1, \mu_2) := \inf \left\{\left(\int_{\X \times \X} \d^2(x,y)\,
d\,q(x,y)\right)^{1/2}:\quad  
q\in \Cpl(\mu_1,\mu_2)\
\right\}
\end{equation}
 on the space $\mathcal P(\X)$ of Borel probability measures on $\X$ where
 $$\Cpl(\mu_1,\mu_2):=\Big\{q\in\mathcal P(\X\times\X), \  (\pi_1)_*q=\mu_1, \ (\pi_2)_*q=\mu_2\Big\}$$
 denotes the set of \emph{couplings} of two probability measures $\mu_1,\mu_2$,
\item the \emph{Boltzmann entropy} 
\begin{equation}
S(\mu)=\Ent (\mu|m) =
\begin{cases}  \int_\X \rho \log
\rho\, d \m, & \mbox{ if} \ \mu=\rho\cdot \m\\
+\,\infty, & \mbox{ if} \ \mu \not\ll \m,
\end{cases}
\end{equation}
regarded as a functional on $\mathcal P(\X)$.
\end{itemize}
%Note that 
The first of these quantities is defined merely using the metric $\d$ on $\X$, the second one merely using the measure $\m$ on $\X$.

\begin{remark}
En passant, we record some nice properties of the underlying metric $\d$ on $\X$ which carry over to the Kantorovich--Wasserstein metric on the \emph{Wasserstein space}
$
{\cal P}_2(X)=\{\mu\in\mathcal P(\X): \ 
\int_\X \d^2(x,x_0)\,\mu(dx) < \infty\}$:
\begin{itemize}
\item $({\cal P}_2(\X), W_2)$ is a complete separable
metric space,

\item $({\cal P}_2(\X), W_2)$ is a
\emph{compact} space or a
\emph{length} space or an
\emph{Alexandrov} space with curvature
$\ge0$ if and only if $(\X,\d)$ is so.
\end{itemize}
\end{remark}

\subsection{Synthetic Ricci Bounds for Metric Measure Spaces}

Following \cite{St06,St06a,LV09}, we now present the so-called \emph{Curvature-Dimension Condition} $\CD(K,N)$, to be considered as a synthetic formulation for ``Ricci curvature $\ge K$ and dimension $\le N$''.
For convenience, we first treat the case $N=\infty$ where 
no constraint on the dimension is imposed.

\begin{definition}
We say that a metric measure space  $(\X,\d,\m)$  has \emph{Ricci curvature $\ge K$} or that it satisfies the
\emph{Curvature-Dimension Condition} $\CD(K,\infty)$
iff
$\forall \mu_0,\mu_1 \in \Pz_2(\X):\ \exists$\
$W_2$-geodesic $(\mu_t)_{t\in[0,1]}$ connecting them s.t. 
\begin{eqnarray}\label{ent-conv}
S(\mu_t)
 \leq  (1-t)\,S(\mu_0) + t\,S(\mu_1) -
\frac{K}{2}\,t(1-t)\, W_2^2(\mu_0,\mu_1).\end{eqnarray}
\end{definition}
\begin{minipage}[c]{0.68\textwidth}
\begin{remark} In other words, the $\CD(K,\infty)$-condition holds true if and only if the Boltzmann entropy is \emph{weakly $K$-convex} on $\Pz_2(\X)$. Recall that $S$ is called \emph{$K$-convex} on $\Pz_2(\X)$ iff \eqref{ent-conv} holds true for \emph{all} $W_2$-geodesics $(\mu_t)_{t\in[0,1]}$ in $\Pz_2(\X)$. The reason for requiring the weaker version is the stability under convergence of the latter (see below).
\end{remark}
\end{minipage}
\begin{minipage}[c]{0.28\textwidth}
\includegraphics[width=4.5cm]{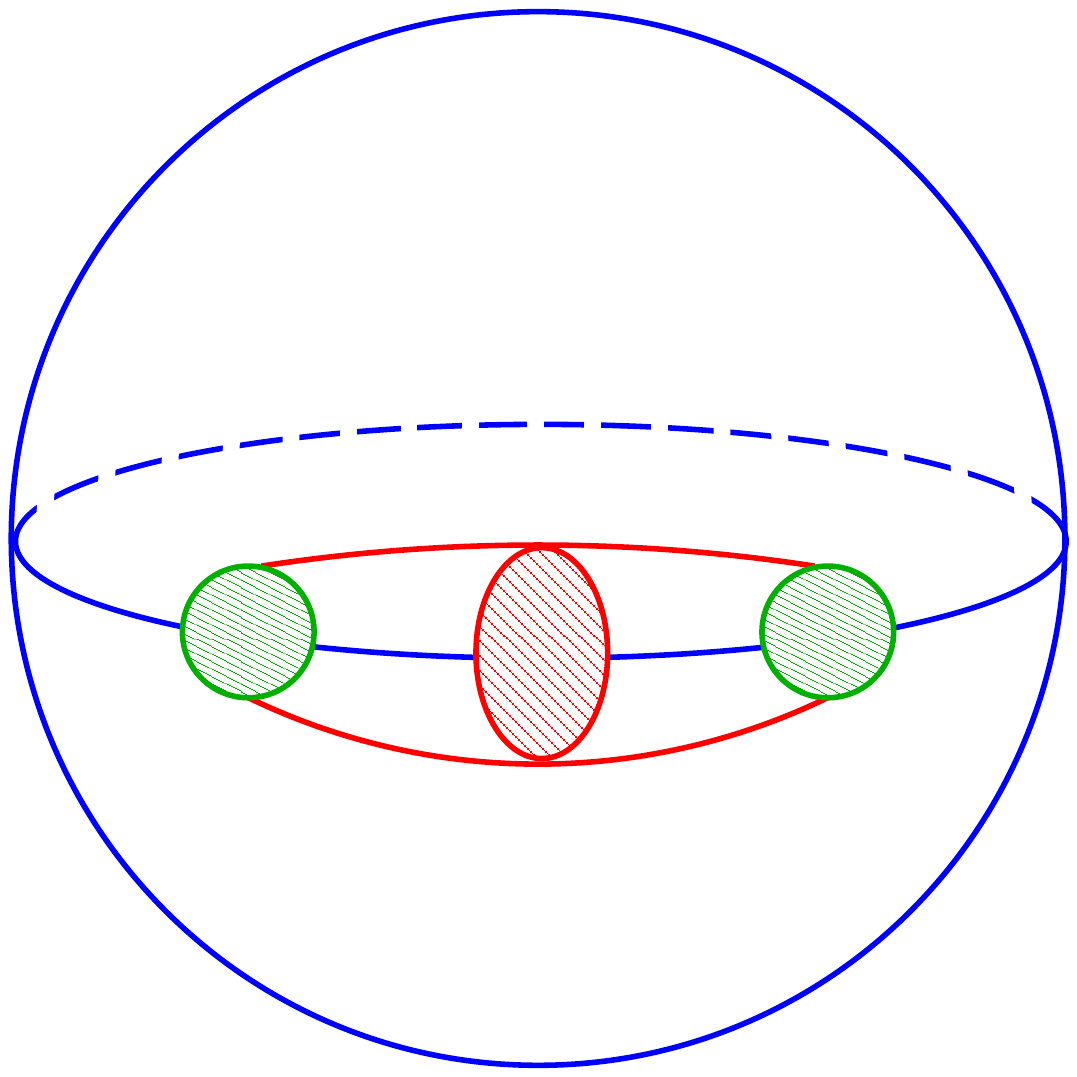}
\end{minipage}

\bigskip

The second case which allows for an easy formulation is $K=0$. Here for finite $N\in\R_+$ the formulation is based on the Renyi type entropy
$$S_N(\nu|\m) :=- \int_\X \rho^{1-1/N}
\, d\m\quad\text{for 
 }\nu=\rho\cdot \m+\nu_s.$$

\begin{definition}
We say that 
$(\X,\d,\m)$ 
 satisfies the  \emph{Curvature-Dimension Condition $\CD(0,N)$} iff
$ \forall \mu_0,\mu_1 \in \Pz_2(\X):\ \exists$
$W_2$-geodesic $(\mu_t)_{t\in[0,1]}$ connecting them  s.t.
\begin{eqnarray} S_N(\mu_t|\m)
& \leq&  (1-t)S_N(\mu_0|\m) + t\,S_N(\mu_1|\m).\end{eqnarray}
\end{definition}

\begin{remark}  It is quite instructive to observe that
 $$S_N(\nu|\m) =-\m(A)^{1/N} \quad
\text{if $\nu$ is unif. distrib. on }A\subset X.$$

\vskip-.6cm
\includegraphics[width=5cm]{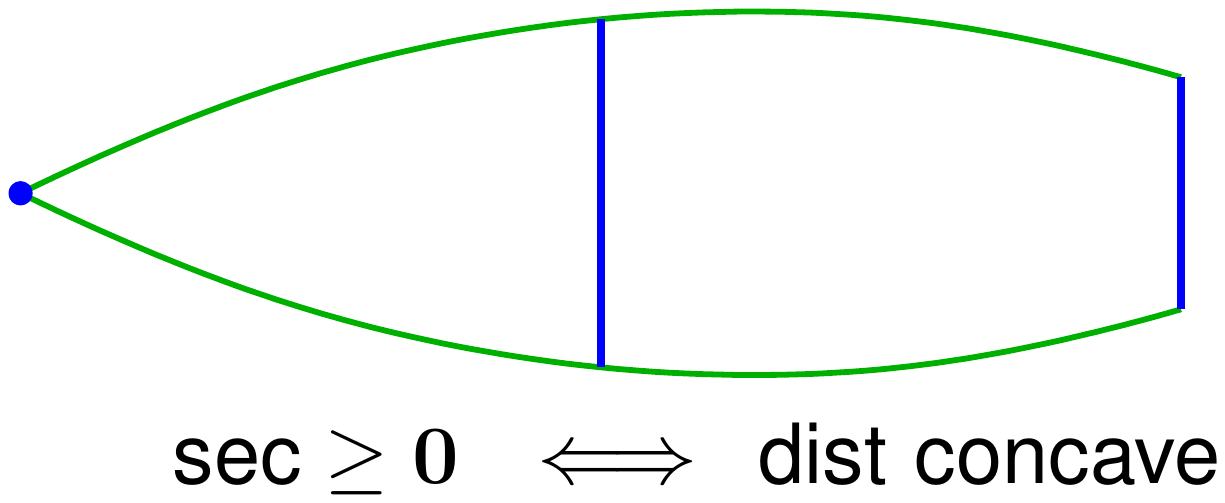}
\includegraphics[width=5cm]{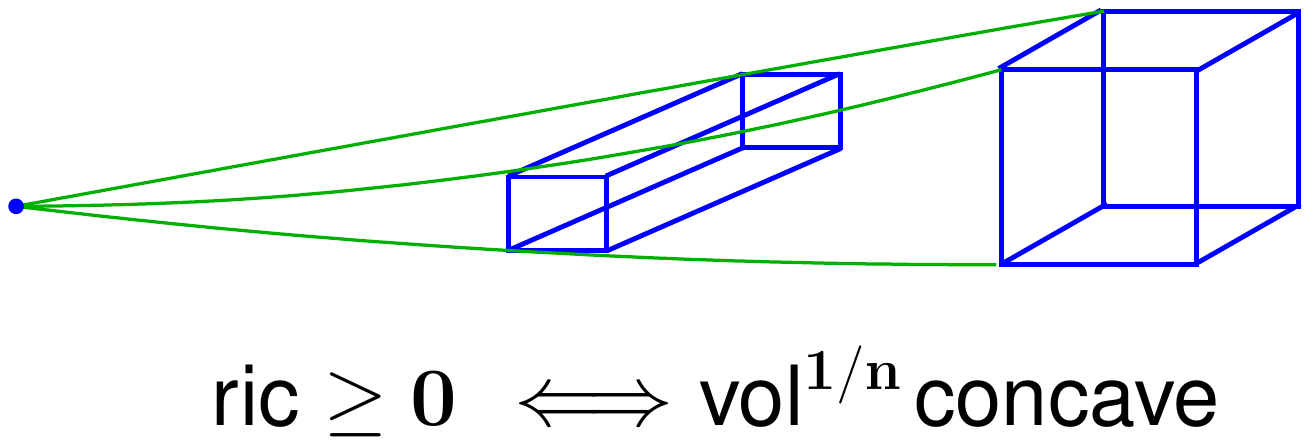}\\[-1.1cm]
Thus the  Curvature-Dimension Condition $\CD(0,N)$ can be vaguely interpreted as a kind of concavity property for the $N$-th root of the volume. This should be seen in context with the facts that i) on $N$-dimensional spaces, the $N$-th root of the volume has the dimension of a length, ii) nonnegative sectional curvature in the sense of Alexandrov can be regarded as a concavity property of distances, and iii) Ricci curvature 
should be regarded as the average of the sectional curvatures. 
\end{remark}

\subsection{The Curvature-Dimension Condition $\CD(K,N)$}

The Curvature-Dimension Condition $\CD(K,N)$ for general pairs of $K,N$ is more involved. It was introduced in \cite{St06a}.
 (Based on that, later on Lott--Villani \cite{LV07} also introduced a slight modification of  it --- the difference, however, will be irrelevant for the sequel. 
In their original paper \cite{LV09}, they   consider only the case $K/N=0$ where the effects of dimension and curvature are decoupled.)

\begin{definition} An mm-space
$(\X,\d,\m)$ satisfies the \emph{Curvature-Dimension Condition}
$\CD(K,N)$ for $K,N \in \mathbb{R}$ (with $N\ge1$) iff
 $ \forall\ \rho_0\m,\ \rho_1\m\in\Pz_2(\X):\ \exists\ W_2\mbox{-geodesic}\ (\rho_t\m)_{t\in[0,1]}$ connecting them and a $W_2$-optimal coupling $q$ of them s.t.
\begin{eqnarray}\label{cd-tau}
\int_\X\rho_t^{1-1/N}(z)\, d\m(z)  & \geq &
\int_{\X\times \X}\left[
{\tau_{K,N}^{(1-t)}(\gamma_0, \gamma_1)}\cdot
\rho_0^{-\,1/N}(\gamma_0)\right.\\
&&\left.\qquad
 + {\tau_{K,N}^{(t)}(\gamma_0, \gamma_1)}\cdot
\rho_1^{-\,1/N}(\gamma_1)\right]\, dq(\gamma_0,\gamma_1).
\nonumber\end{eqnarray}
Here the \emph{distortion coefficients} are given by
$$\tau_{K,N}^{(t)}(x,y):= t^{\frac1N} \left(
\frac{\sin\left(\sqrt{\frac{K}{N-1}}\, t\,\mathrm{d}(x,y)\right)}{\sin\left(\sqrt{\frac{K}{N-1}}\,\mathrm{d}(x,y)\right)}\right)^{\frac{N-1}{N}}$$
in case $K>0$, analogous formula with $\sin\sqrt{\ldots}$ replaced by $\sinh\sqrt{-\ldots}$ in case $K<0$, and $\tau_{K,N}^{(t)}(x,y):=t
$ in case $K=0$.
\end{definition}

The interpretation of $\CD(K,N)$ as a synthetic formulation for  ``Ricci curvature $\ge K$, dimension $\le N$''
is justified by the Riemannian case.

\begin{theorem}[\cite{StvRe05,St05,St06a}]
For Riemannian manifolds $(\M,\g)$,
$$\CD(K,N)\quad \Longleftrightarrow
\quad
\mathrm{Ric}_\M \geq K \quad \mbox{and} \quad \mathrm{dim}_\M \leq N.$$
\end{theorem}

Further examples of metric measure spaces satisfying a $\CD(K,N)$-condition include weighted Riemannian spaces, Ricci limit spaces, Alexandrov spaces, and Finsler spaces.
If one slightly extends the concept of `metric' towards 'pseudo metric', it also includes path spaces (e.g. the
Wiener space with $K=1, N=\infty$) and configuration spaces.

Moreover, many further examples are obtained by \emph{constructions} as
limits, products, cones, suspensions, or warped products.

\section{Geometric Aspects}

The broad interest in --- and the great success of --- 
the concept of the Curvature-Dimension Condition $\CD(K,N)$ is due to  
\begin{itemize}
\item its equivalence to classical lower Ricci bounds in the Riemannian setting, 
\item its stability under convergence and under various constructions, and
\item the fact that it implies almost all of the geometric and functional analytic estimates (with sharp constants!) from Riemannian geometry which depend only on (the dimension and on) lower bounds on the Ricci curvature.
\end{itemize}

\subsection{Volume Growth}

Let us summarize some of the most fundamental geometric estimates.

\begin{theorem}[Bonnet-Myers Diameter Bound  \cite{St06a}] The $\CD(K,N)$-condition with finite $N$ and positive $K$ implies compactness of $\X$ and
\begin{equation}\mathrm{diam}(\X) \quad \leq \quad \sqrt{\frac{N-1}{K}}\cdot
\pi.\end{equation}
\end{theorem}

\begin{theorem}[Bishop-Gromov Volume Growth Estimate  \cite{St06a}] Under $\CD(K,N)$ with finite $N$, for every $x_0\in\X$, the volume growth function $r\mapsto \m(B_r(x_0))$ is absolutely continuous and its weak derivative  $s(r):= \frac\partial{\partial r}\m(B_r(x_0))$ satisfies
 \begin{equation}\label{BG-s}
 {s(r)}\big/ {s(R)} \quad \geq \quad
{\sin\left(\sqrt{\frac{K}{N-1}}r\right)^{N-1}}\Bigg/{\sin\left(\sqrt{\frac{K}{N-1}}R\right)^{N-1}}
\end{equation}
for all $0<r<R$ with the usual re-interpretation of the RHS if $K\le0$ (i.e.~replacing all $\sin(\sqrt{K}\ldots)$ by $\sinh(\sqrt{-K}\ldots)$ in the case $K<0$).
\end{theorem}
As in the smooth Riemannian setting, this differential inequality immediately implies the integrated version:
 \begin{equation*}\label{BG-s}
 \frac{\m(B_r(x_0))}{\m(B_R(x_0))} \quad \geq \quad
\frac{\int_0^r{\sin\left(\sqrt{\frac{K}{N-1}}t\right)^{N-1}}dt}{\int_0^R{\sin\left(\sqrt{\frac{K}{N-1}}t\right)^{N-1}}dt}
\end{equation*}
for all $0<r<R$, and thus in particular
$$\m(B_R(x_0))\le C\,r^N\,\exp\big(\sqrt{(N-1)K^-} \,R\big).$$
The results so far assumed that $N$ is finite.
In the case $N=\infty$, the $\CD(K,N)$-condition implies a \emph{novel volume growth estimate} \cite{St06}, not known before in the Riemannian setting,
\begin{equation}\m(B_R(x_0))\le \exp\bigg(\frac {K^-}2 R^2+c_1R+c_0\bigg).
\end{equation}
It can be seen as complementary to the concentration of measure phenomenon. The sharpness is illustrated by the following
\begin{example} Consider $\X=\R, \d=|\,.\,|$ and $d\m(x)=\exp(\frac\kappa2|x|^2)$ for %some 
$\kappa>0$. Then $(\X,\d,\m)$ satisfies $\CD(-\kappa,\infty)$, and $\m(B_R(x))\ge \exp(\frac\kappa2(R-\frac12)^2)$ for all $x$ and $R\ge \frac12$.
\end{example}
The Curvature-Dimension Condition $\CD(K,N)$ also implies numerous further geometric estimates, among them
the \emph{Brunn-Minkowski Inequality} \cite{St06a} and the  \emph{Borell-Brascamp-Lieb Inequality} \cite{Bacher10}.
What remained an open problem for many years was the L\'evy--Gromov Isoperimetric Inequality which only recently was proven by Cavalletti--Mondino.

\begin{theorem}[L\'evy-Gromov Isoperimetric Inequality \cite{CavallettiMondino17}]  Let $(\X,\d,\m)$ be an mm-space which satisfies $\CD(K,N)$ and let $\hat\X$ be a $\CD(K,N)$-model space. Then for every subset $E\subset \X$ and every  spherical cap $B\subset\hat \X$,
\begin{equation}\frac{|\partial E|}{|\X|}\ge \frac{|\partial B|}{ |\hat \X|}\quad\text{if}\quad 
\frac{|E|}{|\X|}=\frac{|B|}{|\hat \X|}.
\end{equation}
Here $|\,.\,|$ denotes the respective volume or surface measure.
\end{theorem}

\subsection{The  Space of Spaces}

Two mm$_1$-spaces will be called \emph{isomorphic} --- and henceforth identified --- iff there exists a measure preserving isometry between the supports of the respective measures.
It is a quite remarkable observation that the space 
$\boldsymbol{\rm\Xi}$ of isomorphism classes of normalized mm$_1$-spaces itself is a geodesic space.

The \emph{$L^{p}$-transportation distance} between two mm$_1$-spaces $(\X_0,\d_0,\m_0)$ and $(\X_1,\d_1,\m_1)$  is
	defined  for $p\in[1,\infty)$ as
		$${\mathbb D}_p\Big((\X_0,\d_0,\m_0),(\X_1,\d_1,\m_1)\Big)
= \inf_{\d,\m}   \bigg( \int_{\X_0\times \X_1} 
		 \d(x_0,x_1)^p d{\m}
		(x_0,x_1)\bigg)^{1/p}$$
	where the infimum is taken over all {\emph{couplings $\m$} of $\m_0$ and $\m_1$ 
and over all \emph{couplings $\d$} of $\d_0$ and $\d_1$ (i.e. metrics on $\X_0\sqcup \X_1$ which coincide with $\d_0$ on $\X_0$ and with $\d_1$ on $\X_1$),  \cite{St06}.
With slight modifications, this definition also extends to 
 $p=\infty, p\in (0,1)$ and $p=0$.

A closely related concept is the \emph{$L^p$-distortion  distance} between  mm$_1$-spaces 
	defined  for $p\in[1,\infty)$ as
	\begin{eqnarray*}
	\lefteqn{\boldsymbol{\rm\Delta}_p\Big((\X_0,\d_0,\m_0),(\X_1,\d_1,\m_1)\Big)}\\
&=& \inf_{\m}   \bigg( \int_{\X_0\times \X_1} \int_{\X_0\times \X_1}
		\Big| \d_0(x_0,y_0)-\d_1(x_1,y_1) \Big|^p d{\m}
		(x_0,x_1)d{\m}(y_0,y_1)\bigg)^{1/p}
	\end{eqnarray*}
where the infimum is taken over all \emph{couplings} $\m$ of $\m_0$ and $\m_1$, and again with slight modifications also extended to 
 $p=\infty, p\in (0,1)$ and $p=0$.
	Under uniform control of the moments of the involved metric measure spaces, the topologies induced by all these metrics are the same and coincide with that of  \emph{Gromov's box distance} $\Box_\lambda$  and to that of 	measured Gromov-Hausdorff convergence.
\begin{lemma}[\cite{gpw,Memoli11,StSoS12}]
\begin{enumerate}
\item $\forall p\in [0,\infty)$: \ 
${\mathbb D}_p$ is complete whereas $\boldsymbol{\rm\Delta}_p$ is not complete
\item ${\mathbb D}_p$-convergence $\quad\Longleftrightarrow\quad {\mathbb D}_0$-convergence and convergence of $p^{th}$ moments
\item $\boldsymbol{\rm\Delta}_p$-convergence $\quad\Longleftrightarrow\quad \boldsymbol{\rm\Delta}_0$-convergence and convergence of $p^{th}$ moments
\item
${\mathbb D}_0$-convergence $\quad\Longleftrightarrow\quad$ ${\boldsymbol{\rm\Delta}}_0$-convergence $\quad\Longleftrightarrow\quad$ $\Box_\lambda$-convergence.
\end{enumerate}
\end{lemma}

The main result here is that the \emph{space of spaces} is an Alexandrov space. 

\begin{theorem}[\cite{StSoS12,StSoS20}] The metric space $(\boldsymbol{\rm\Xi}_2,\boldsymbol{\rm\Delta}_2)$ of isomorphism classes of mm$_1$-spaces  
 is a \emph{geodesic space} with \emph{nonnegative curvature}.
\end{theorem}

The tangent space  (for the space of spaces) at a given mm$_1$-space
admits an explicit representation and so does the \emph{symmetry group}, the latter e.g.~in terms of \emph{optimal self-couplings}.
Of particular interest are finite dimensional subspaces of the space of spaces.

\begin{proposition} For each $n\in\N$, the subspace of $n$-point spaces (i.e. mm$_1$-spaces with equal mass on $n$-points) is a \emph{Riemannian orbifold} with nonnegative curvature.
\end{proposition}

\subsection{Stability, Compactness}

Converging sequences of mm$_1$-spaces can always be embedded into common metric spaces. The stability of the $\CD(K,N)$-condition then simply amounts to the lower semicontinuity of the Renyi type entropy for weakly convergent sequences of probability measures.

\begin{theorem} 
The curvature-dimension condition is
\emph{stable} under $\mathbb D_0$-convergence of mm$_1$-spaces.
\end{theorem}

The volume growth estimates entailed by the $\CD(K,N)$-condition  together with  the stability of the latter under convergence, allow us to turn Gromov's pre-compactness theorem under Ricci bounds into a compactness theorem.

\begin{theorem}   For every triple $K, N, L\in\mathbb{R}$, the space  of all mm$_1$-spaces
 $(\X,\d,\m)$ that satisfy  $\CD(K,N)$ and have diameter $\leq
L$ is \emph{compact}.
\end{theorem}

\subsection{Local to Global}

A crucial property of curvature bounds both in Riemannian geometry and in the geometry of Alexandrov spaces is the \emph{local-to-global property}: sharp global estimates follow from uniform local curvature assumptions.
For the synthetic Ricci bounds for mm-spaces, this 
is a highly non-trivial claim.
To deal with it, we restrict ourselves to non-branching geodesic spaces.

The first  \emph{Globalization Theorem} was obtained in the case $K/N=0$ where curvature and dimension effects are de-coupled.
\begin{proposition}[\cite{St06,St06a,LV09}] \ If $K=0$ or $N=\infty$ then every mm-space $(\X,\d,\m)$  satisfies 
$$\CD(K,N) \text{  locally }
 \quad \Longleftrightarrow \quad \CD(K,N) \text{ globally}.$$ 
\end{proposition}

Further progress then was based on the   \emph{Reduced Curvature-Dimension Condition} $\CD^*(K,N)$ defined similarly as $\CD(K,N)$ but now with the distortion coefficient $\tau_{K,N}^{(t)}(x,y)$ in \eqref{cd-tau} replaced by the reduced coefficients
$\sigma_{K,N}^{(t)}(x,y):=
{\sin\left(\sqrt{\frac{K}{N}}t\,\mathrm{d}(x,y)\right)}\Big/{\sin\left(\sqrt{\frac{K}{N}}\mathrm{d}(x,y)\right)}.$

\begin{proposition}[\cite{BacherSt10}] For all $K,N\in\R$ and all mm-spaces,
 $$\CD(K,N) \text{  locally }  \quad \Longleftrightarrow \quad  \CD^*(K,N) \text{  locally }
 \quad \Longleftrightarrow \quad \CD^*(K,N) \text{ globally}.$$ 
\end{proposition}

Only recently, the Globalization Theorem could be proven in full generality by Cavalletti--Milman  (with a minor extension by Zhenhao Li  removing the finiteness assumption for the underlying measure). Their approach is based on Klartag's \cite{Klartag17} \emph{needle decomposition} and  the \emph{localization technique} developed by Cavalletti--Mondino \cite{CavallettiMondino17}.
\begin{theorem}[\cite{CavallettiMilman21}, {[Li21+]}] $$\CD(K,N) \text{ locally}\quad\Longleftrightarrow\quad \CD(K,N) \text{ globally.}$$
\end{theorem}

\section{Analytic Aspects}
A deeper understanding of the role of synthetic lower Ricci bounds on singular spaces will be obtained through links with 
spectral properties of the Laplacians and estimates for heat kernels on such spaces

\subsection{Heat Flow on Metric Measure  Spaces}

There are two different (seemingly unrelated) approaches to define  the  \emph{heat equation} on an mm-space $(\X,\d,\m)$
\begin{itemize}
\item[$\ast$]
either as gradient flow in $L^2(\X,\m)$ for the
\emph{energy}
$${\mathcal{E}(u)=\frac12\int_X |\nabla u|^2\,dm=
\liminf_{v \to u \mathrm{\, in \, }L^2} \frac12\int_X (\mathrm{lip}_x v)^2\,dm(x)}
$$
with $\mathrm{lip}_x v(x)=\limsup\limits_{y\to x}\frac{|v(x)-v(y)|}{\d(x,y)}$ 
and
$|\nabla u|=$  minimal weak upper gradient,
\item[$\ast$]
or  as gradient flow in
$\mathcal{P}_2(\X)$ for the
\emph{Boltzmann entropy}
$$\mbox{Ent}(u)=\int_X u\log u\,dm.$$
\end{itemize}
The former approach (the traditional point of view) has the advantage that the energy --- if it exists ---  is always convex and thus guarantees the existence of the gradient flow. Its disadvantage is that it relies on the concept of weakly differentiable functions. However, all analytic problems related to the notion of energy have been fully resolved in the trilogy \cite{AGS1,AGS2,AGS3}  by Ambrosio--Gigli--Savar\'e.

The latter approach (the novel perspective of Otto) has the advantage that the entropy is always obviously  well-defined. However, for its gradient flow to exist, additional assumptions are required, e.g.~that the  entropy is semi-convex. Up to minor technicalities, this simply says that the underlying mm-space has lower bounded synthetic Ricci curvature. Under this minimal assumption, indeed, 
both approaches coincide.

\begin{theorem}[AGS1] For every mm-space  $(\X,\d,\m)$  that satisfies $\CD(K,\infty)$ for some $K\in\R$, the energy approach and the entropy approach coincide.
\end{theorem}

\begin{example} There are plenty of examples to which this result applies. 
The most prominent among them (and the authors who first proved it) are
\begin{enumerate}
\item[i)] \emph{Euclidean space}  $\mathbb{R}^n$:  \ Jordan--Kinderlehrer--Otto \cite{JKO98}
\item[ii)] \emph{Riemann manifolds} $(\M,\g)$:  \ Ohta \cite{Ohta09}, Savar\'e \cite{Savare07}, Villani \cite{VillaniBook2}, Erbar \cite{Erbar10}

\item[iii)] \emph{Finsler spaces $(\M,\mathsf F,\m)$}: \ Ohta--Sturm \cite{OhtaSt09}
\item[iv)] \emph{Alexandrov spaces}: \ Gigli--Kuwada--Ohta \cite{GKO13}
\end{enumerate}
\end{example}

\begin{example} In many other cases not covered by any $\CD$-condition we know that the energy approach and the entropy approach coincide:

\begin{enumerate}
\item \emph{Heisenberg group} (unbounded curvature): \ Juillet \cite{Juillet14}
\item \emph{Wiener space} (degenerate distance): \ Fang--Shao--Sturm \cite{FangShaoSturm}
\item \emph{Configuration space} (degenerate distance): \ Erbar--Huesmann \cite{ErbarHuesmann15}
\item \emph{Neumann Laplacian} (unbounded curvature if nonconvex): \ Lierl--Sturm \cite{LierlSturm}
\item \emph{Dirichlet Laplacian} (no mass conservation): \ Profeta--Sturm \cite{ProfetaSturm}

\item \emph{Discrete spaces} (no $W_2$-geodesics): \ Maas \cite{Maas11}, Mielke \cite{Mielke11}
\item \emph{Levy semigroups} (no $W_2$-geodesics): \ Erbar \cite{ErbarLevy}
\item \emph{Metric graphs} (unbounded curvature): \ Erbar--Forkert--Maas--Mugnolo \cite{EFMM21}.
\end{enumerate}
In the latter examples e), f), and g), the concept of ``gradient flow for the Boltzmann entropy'' has to be slightly adapted.
\end{example}

\subsection{Curvature-Dimension Condition: Eulerian vs. Lagrangian}

 Besides the Lagrangian formulation of synthetic Ricci bounds in terms of semiconvexity properties of the entropy, there is also an Eulerian formulation in terms of the energy: the celebrated \emph{curvature-dimension (or $\Gamma_2$) condition of Bakry--\'Emery}. 
 It is a groundbreaking observation that both formulations  
 are equivalent in great generality.

 For this equivalence to hold, we now make the standing assumption that $(\X,\d,\m)$ is \emph{infinitesimally Hilbertian}, 
i.e.~the energy $\mathcal E$ is quadratic  or, in other words, Laplacian and  heat flow are linear. 
For convenience we will also assume that the mm-space under consideration has the \emph{Sobolev-to-Lipschitz property} and volume growth bounded by $\ee^{Cr^2}$. Note that both of these latter properties follow from the validity of the Lagrangian $\CD(K,N)$-condition.

\begin{theorem}[\cite{AGS2,AGS3,EKS}] \label{4equiv}
Under the above assumptions, the following properties are equivalent:
\begin{enumerate}
\item[\rm (i)] the synthetic Ricci bound $\CD(K,N)$, briefly reformulated as
 $$\mbox{\rm Hess}\ S -\frac1N (\nabla S)^{\otimes 2}\ge K\quad\text{on }(\Pz_2(\X), W_2),$$
 \item[\rm (ii)] the transport estimate
 $$W_2^2(P_s\mu,P_t\nu)\le e^{-K\tau}\, W_2^2(\mu,\nu) +2N\frac{1-e^{-K\tau}}{K\tau}(\sqrt s-\sqrt t)^2
$$
with $\tau:=\frac23(s+\sqrt{st}+t)$,
\item[\rm (iii)] the gradient estimate
$$ \left|\nabla P_tu\right|^2+\frac{4Kt^2}{N(e^{2Kt}-1)}
\left|\Delta P_tu\right|^2
\le  
e^{-2Kt}  P_t\left|\nabla u\right|^2,
$$
\item[\rm (iv)] the Bochner inequality
$$\frac12\Delta|\nabla u|^2-\langle\nabla u, \nabla\Delta u\rangle\ge K \cdot |\nabla u|^2 \, +\frac1N(\Delta u)^2,$$
also known as Bakry--\'Emery criterion and written in comprehensive form as 
$$\Gamma_2(u)\ge K\cdot \Gamma(u)+\frac1N(\Delta u)^2.$$
\end{enumerate}
\end{theorem}
 These equivalences allow for easy explanations and/or intuitive interpretations. 
The equivalence (iii) $\Leftrightarrow$ (iv), indeed, is known since decades as a basic result of the so-called  $\Gamma$-calculus of Markov semigroups \cite{Bakry85,Bakry94},
and easily follows by differentiating $s\mapsto P_{t-s}(|\nabla P_su|^2)$.
 The equivalence (i) $\Leftrightarrow$ (ii), from an heuristic point of view, is a consequence of the fact that the heat flow is the gradient flow for the entropy w.r.t.~the metric $W_2$.
 Finally, the equivalence (ii) $\Leftrightarrow$ (iii) is the important \emph{Kuwada duality} which extends the celebrated \emph{Kantorovich--Rubinstein duality} towards $p\not=1, q\not=\infty$.
 
The rigorous proofs of  the above equivalences by Ambrosio--Gigli--Savar\'e \cite{AGS2,AGS3} (for the case $N=\infty$) and Erbar--Kuwada--Sturm \cite{EKS} (for the general case) are rather sophisticated and mark  crucial milestones in the development of the theory. 

\begin{remark} The Bakry--\'Emery estimate $$\Gamma_2(u)- K\cdot |\nabla u|^2\ge\frac1N(\Delta u)^2\qquad(\forall u)$$ has a remarkable \emph{self-improvement property} \cite{Bakry85,Bakry94,BakryQian,SavareSelf,ErbarSturm22} asserting that it implies the seemingly stronger estimate 
\begin{align*}\Gamma_2(u)- K\cdot |\nabla u|^2&\ge\frac1N(\Delta u)^2+\frac N{N-1}\Big| \big|\nabla|\nabla u|\big|-\frac1N|\Delta u|\Big|^2\\
&= \big|\nabla|\nabla u|\big|+\frac1{N-1}\Big| \big|\nabla|\nabla u|\big|-|\Delta u|\Big|^2
\qquad\qquad(\forall u).
\end{align*}
This leads to improved gradient estimates and improved transport estimates which e.g.~in the case  $N=\infty$ read as
$$ \left|\nabla P_tu\right|
\le  
e^{-Kt}  P_t\left|\nabla u\right|,\qquad
W_\infty(P_t\mu,P_t\nu)\le e^{-K t}\, W_\infty(\mu,\nu).
$$

\end{remark}

\subsection{$\RCD(K,N)$-Spaces \quad --- \quad Functional Inequalities}

We will say that an mm-space satisfies \emph{the $\RCD(K,N)$-condition} iff it satisfies the $\CD(K,N)$-condition and iff it is infinitesimally Hilbertian. For these mm-spaces, the full machinery of geometric analysis and Riemannian calculus can be developed and  far reaching structural assertions can be derived.

Here we have to restrict ourselves to present only a selection of the many results proven so far. 
And  we will not formulate detailed estimates (except for the first result), we will just mention the respective results.

\begin{theorem} The following estimates hold true (each of them with sharp constants) on any mm-space which satisfies an $\RCD(K,N)$-condition for some $K\in\R$ and for
 $N\le\infty$:
\begin{itemize}
\item[$\ast$] Poincar\'e/Lichnerowicz inequality  \cite{LV07}: \quad  $\lambda_1\ge \frac N{N-1}\,K,$
\end{itemize}
moreover, for $N<\infty$:
\begin{itemize}
\item[$\ast$]  Laplace comparison \cite{Gigli15differential},
\item[$\ast$] Bochner inequality  \cite{EKS},
\item[$\ast$] Li-Yau differential Harnack inequality,  Gaussian heat kernel estimates \cite{GarofaloMondino},
\item[$\ast$] Sobolev, Cheeger, and Buser inequalities  \cite{Profeta15}, \cite{DePontiMondino},
\end{itemize}
whereas   for $N=\infty$:
\begin{itemize}
\item[$\ast$]     Talagrand  and Logarithmic Sobolev inequalities \cite{LV07},
\item[$\ast$]  Wang Harnack inequality \cite{LiHuaiqian}, upper Gaussian heat kernel estimate \cite{Tamanini2019harnack},
and
Ledoux inequality \cite{DePontiMondino}.
\end{itemize}
 \end{theorem}

In all the previous results, the dimensional parameter has always been a number  $N\ge1$ (which in turn then even implies that 
$N\ge \text{dim}_{\mathcal H}(\X)$). Quite remarkably,  various of these results also admit versions 
where the  \emph{dimensional parameter $N$ is a negative number}, see e.g.~\cite{ohta2016k,ohta2011displacement,milman2017beyond,magnabosco2021optimal,magnabosco2021convergence}.

\subsection{$\RCD(K,N)$-Spaces \quad --- \quad Splitting and Rigidity}

In the smooth Riemannian setting, an important consequence of nonnegative Ricci curvature is the Cheeger-Gromoll splitting theorem.
In order to extend this  to metric measure spaces, it is essential to assume that the underlying spaces are infinitesimally Hilbertian.

\begin{theorem}[Splitting Theorem \cite{Gigli14Splitting}]
If an mm-space $(\X,\d,\m)$ satisfies $\RCD(0,N)$ and \emph{contains a line} then 
$\X={\mathbb R}\times \X'$ for some $\RCD(0,N-1)$-space $(\X',\d',\m')$.
\end{theorem}

The counterpart to the splitting theorem for positive lower Ricci bound is Cheng's maximal diameter theorem.

\begin{theorem}[Maximal Diameter  Theorem \cite{ketterer2015cones}]
If an mm-space $(\X,\d,\m)$ satisfies $\RCD(N-1,N)$ and has \emph{diameter $\pi$} then $\X$ is the spherical suspension of 
some $\RCD(N-2,N-1)$-space $(\X',\d',\m')$.
\end{theorem}
In the smooth Riemannian setting,  the maximal diameter theorem provides a more far-reaching conclusion, namely, that $\X$ is the round $N$-sphere. In the singular setting, however, this conclusion is false \cite{ketterer2015cones}.

On the other hand, such a far-reaching conclusion can be drawn from the maximality of the spherical size.
\begin{theorem}[Maximal Spherical Size Theorem \cite{ErbarSturm20}] If an mm-space $(\X,\d,\m)$ satisfies  $\RCD(N-1,N)$ and
\begin{equation}-\int_\X\int_\X\cos\Big( \d(x,y)\Big)\,d\m(x)\,d\m(y)\ge0,\end{equation}
then $N\in\mathbb N$ and
$(\X,\d,\m)$ is isomorphic to the $N$-dimensional round sphere ${\mathbb S}^N$.
\end{theorem}

Closely related to the maximal diameter theorem is Obata's theorem on the minimality of the spectral gap. 
\begin{theorem}[Obata's Theorem \cite{ketterer2015obata}]
If an mm-space $(\X,\d,\m)$ satisfies $\RCD(N-1,N)$ and has \emph{spectral gap $N$} then $\X$ is the spherical suspension of 
some $\RCD(N-2,N-1)$-space $(\X',\d',\m')$.
\end{theorem}

This splitting theorem indeed also admits an extension to $N=\infty$ which states that an 
mm-space $(\X,\d,\m)$ that satisfies $\RCD(1,\infty)$ and has spectral gap $1$ splits off a Gaussian factor \cite{gigli2020rigidity}.

\subsection{$\RCD(K,N)$-Spaces \quad --- \quad Structure Theory}

Since blow-ups of $\RCD(K,N)$-spaces are $\RCD(0,N)$-spaces which contain lines, a 
sophisticated 
iterated application of the splitting theorem will
lead to deep insights into tangent spaces and local structure of $\RCD$-spaces.

\begin{theorem} [Rectifiability and Constancy of Dimension \cite{mondino2019structure,brue2020constancy}]
 If $(\X,\d,\m)$ satisfies $\RCD(K,N)$ then  
\begin{enumerate}
\item
$\X=\bigcup_{k=1}^{\lfloor N\rfloor} {\mathcal R}_k \cup {\mathcal N}, \quad m(\mathcal N)=0$,
\item each  $\mathcal R_k$ is covered by countably many measurable sets which are $(1+\epsilon)$-biLipschitz equivalent to subsets of $\R^k$,
\item $\m$ and $\mathcal H^k$ are mutually abs.~cont.~on $\mathcal R_k$,
\end{enumerate}
and even more,
\begin{enumerate}
\item[\rm d)] $\exists n\in{\mathbb N}$ such that $m(\mathcal R_k)=0$ for all $k\not=n$.
\end{enumerate}
\end{theorem}

Besides the two landmark contributions to this structure theory mentioned above, numerous 
important results were obtained \cite{kell2016volume,gigli2016behaviour,de2017conjecture,ambrosio2018short}.
Particularly nice insights could be obtained in the case $N=2$.

%Bru\`e/Pasqualetto/Semola '20, 

\begin{corollary}[\cite{lytchak2018ricci}]%, \ Kapovich/Ketterer '19) }
$\RCD(K,2)$-spaces with $\m=\mathcal H^2$ are  Alexandrov spaces.
\end{corollary}
 
Further challenges then concern the boundaries of mm-spaces. Various concepts how to define them and related results were presented in \cite{de2018non,kapovitch2021topology,kapovitch2020metric}.
Important contributions to the analysis of tangent cones and to the regularity theory for non-collapsed RCD-spaces were provided in
\cite{kitabeppu2019sufficient,antonelli2019volume,honda2020new}.
Based on these results, a precise description could be derived.

\begin{theorem}[\cite{brue2020boundary}]
Let $(\X,\d,\m)$ be a non-collapsed $\RCD(K,N)$-space (with $\m=\mathcal H^N, N\in\mathbb N$). Then
\begin{enumerate}
\item $\exists$ stratification $\mathcal S^0\subset  \mathcal S^0\subset  \ldots \subset \mathcal S^{N-1}=\mathcal S=\X\setminus\mathcal R_N$,
\item the boundary $\partial \X:=\overline{\mathcal S^{N-1}\setminus \mathcal S^{N-2}}$ is $(N-1)$-rectifiable,
\item $T_x\X\simeq \R^{N-1}\times\R_+$ for $x\in \mathcal S^{N-1}\setminus \mathcal S^{N-2}$,
\item $\X\setminus \mathcal S^{N-2}$ is a topological manifold with boundary. 
\end{enumerate}
 \end{theorem}

\section{Recent Developments}
The concept of synthetic Ricci bounds for singular spaces
turned out to be extremely fruitful, both for theory and applications. A rich theory of mm-spaces satisfying such uniform lower Ricci bounds has been established.
The last 15 years have seen a wave of impressive results --- many of them going far beyond the previously described scope.

In the following, we will first present in detail 
recent developments concerning
\begin{itemize}
\item[$\ast$] heat flow on time-dependent mm-spaces and super-Ricci flows,
\item[$\ast$] 
second order calculus, upper Ricci bounds, and transformation formulas,
\item[$\ast$] distribution-valued lower Ricci bounds, 
\end{itemize}
and then 
briefly summarize several further developments.

\subsection{Heat Flow on Time-dependent MM-Spaces and Super-Ricci Flows}

Whereas construction and properties of the heat flow on `static' mm-spaces $(\X,\d,\m)$ --- in particular, its relation to synthetic lower bounds on the Ricci curvature ---  by now  are well understood in great generality,
analogous questions 
for time-dependent families  of mm-spaces $(\X_t,\d_t,\m_t), t\in I=(0,T)$, until recently remained widely open:
\begin{itemize}
\item How to define a heat propagator $(P_{t,s})_{t\ge s}$ acting on functions in $L^2(\X_s,\m_s)$ and/or its dual $(\hat P_{t,s})_{s\le t}$ acting on measures  on $\X_t$? 

Can they be regarded as gradient flows of (time-dependent) energy or entropy functionals in function/measure spaces with time-dependent norms or metrics?
\item
Is there a parabolic analogue to synthetic lower Ricci bounds? Can one formulate it as ``dynamic convexity'' of a  time-dependent  entropy functional? How is this related to the notion of super-Ricci flows for families of Riemannian manifolds?
\item Are there ``parabolic versions'' of the  functional inequalities that characterize synthetic lower Ricci bounds?
\end{itemize}
Within recent years,  for families of  mm-spaces $(\X,\d_t,\m_t)$, $t\in (0,T)$, such that
\begin{itemize}
\item[$\ast$] for every $t\in I$
 the mm-space $(\X,\d_t,\m_t)$  satisfies an $\RCD(K,N)$-condition,
\item[$\ast$] some regular $t$-dependence of $\d_t$ and $\m_t$,
\end{itemize}
these questions found affirmative answers.

\begin{definition}[\cite{sturm2018super}] A family of mm-spaces
 $(\X,\d_t,\m_t)_{t\in(0,T)}$ is called \emph{{super-Ricci flow}} iff the 
  function $$\Ent: \ (0,T)\times\mathcal P(X)\to(-\infty,\infty], \quad
(t,\mu)\mapsto \Ent_t(\mu):= \Ent(\mu|m_t)$$
 is \emph{dynamically convex}  on $\mathcal P(X)$ --- equipped with the 1-parameter family of metrics $W_t$
(= $L^2$-Kantorovich--Wasserstein metrics w.r.t. $\d_t$) --- 
 in the following sense: for  all $\mu^0,\mu^1$ and a.e.  $t$  there exists a $W_t$-geodesic $(\mu^a)_{a\in[0,1]}$ s.t. 
\begin{equation}
\partial_a \Ent_t(\mu^0)- \partial_a \Ent_t(\mu^1)\le \frac12\partial_t W_t^2(\mu^0,\mu^1).
\end{equation}
\end{definition}
\begin{example}
A family of Riemannian manifolds $(\M,\g_t), t\in (0,T)$ is a 
super-Ricci flow in the previous sense iff $$\mathrm{Ric}_t+\frac12\partial_t\g_t\ge0.$$
 Recall that  $(\M,\g_t)_{t\in (0,T)}$ is called \emph{Ricci flow} if
$\mathrm{Ric}_t+\frac12 \partial_t\g_t=0$.
These properties can be regarded as the 
parabolic analogue to
nonnegative (or vanishing, resp.) Ricci curvature for static manifolds. 
\end{example}

Whereas in the static setting the gradient flow for the energy and the gradient flow for the entropy characterize the same evolution (either in terms of densities or in terms of measures), this is no longer the case in the dynamic setting: here one is characterizing the forward evolution whereas the other one is characterizing the backward evolution.

\begin{theorem}[\cite{KopferSturm18}] In the previous setting, there exists a well-defined  heat propagator $(P_{t,s})_{t\ge s}$ acting on functions in $L^2(\X,\m_s)$ and its dual $(\hat P_{t,s})_{s\le t}$ acting on measures  on $\X$. Moreover,
\begin{enumerate}
\item[$\ast$]
$\forall u\in{\mathcal Dom}(\mathcal E), \forall s\in I$
the heat flow $t\mapsto u_t=P_{t,s} u$ is the unique \emph{forward gradient flow} for the Cheeger energy $\frac12\mathcal E_s$ in $L^2(\X,\m_s)$.\vskip.35cm
\item[$\ast$]
$\forall\mu\in{\mathcal Dom}(\Ent), \forall t\in I$
the dual heat flow $s\mapsto \mu_s=\hat P_{t,s}\mu$ is the unique \emph{backward gradient flow} for the Boltzmann entropy $\Ent_t$ in $({\mathcal P}(\X), W_t)$ 
provided $(\X,\d_t,\m_t)$ is a super-Ricci flow.
\end{enumerate}
\end{theorem}
Both gradient flows can be obtained as limits of corresponding steepest-descend schemes (aka JKO-schemes) adapted to the time-dependent setting \cite{Kopfer}.

\medskip

In 
analogy to Theorem \ref{4equiv}, the Lagrangian characterization of super-Ricci flows (in terms of dynamic convexity of the entropy) 
turns out to be equivalent to an
Eulerian characterization (in terms of a dynamic $\Gamma_2$-inequality), to a gradient estimate for the forward evolution, and to a transport estimate (as well as to a  pathwise Brownian coupling property) for the backward evolution. 

\begin{theorem}[\cite{KopferSturm18}] The following are equivalent
\begin{enumerate} 
\item
$\partial_a \Ent_t(\mu^a)\big|_{a=0}- \partial_a \Ent_t(\mu^a)\big|_{a=1}\le \frac12\partial_t W_t^2(\mu^0,\mu^1)$ \vskip.35cm
 \item
$W_s(\hat P_{t,s}\mu,\hat P_{t,s}\nu)\le W_t(\mu,\nu)$ \vskip.35cm
\item
$\forall x,y, \forall t: \ \exists$  coupled  backward Brownian motions $(X_s,Y_s)_{s\le t}$ starting at $t$ in $(x,y)$ s.t.
$\d_s(X_s,Y_s)\le \d_t(x,y)$
 a.s. for all $s\le t$\vskip.35cm
\item
$|\nabla_t(P_{t,s} u)|^2\le P_{t,s}( |\nabla_su|^2)$ \vskip.35cm
\item
$\Gamma_{2,t}\ge \frac12\partial_t\Gamma_t$
\qquad\qquad\hfill  where {$\Gamma_{2,t}(u)=\frac12\Delta_t|\nabla_t u|^2-\langle \nabla_t u, \nabla_t\Delta_t u\rangle$}.
\end{enumerate}
\end{theorem}

This result in particular extends a previous characterization of super-Ricci flows of \emph{smooth} families of Riemannian manifolds in terms of the previous assertion b) by McCann--Topping \cite{mccann2010ricci} and in terms of the previous assertion c) by Arnaudon--Coulibaly--Thalmaier \cite{arnaudon2008brownian}.

There is a whole zoo of further functional inequalities which characterize super-Ricci flows. Several implications for  the subsequent assertions  were new even in the static case. %For well-known similar equivalences in the smooth case, see e.g.~Wang.

\begin{theorem}[\cite{KopferSturm21}]
Each of the following assertions is equivalent to  any of the above or, in other words, to $(\X,\d_t,\m_t)_{t\in I}$ being a  super-Ricci flow
\begin{enumerate} 
\item[\rm  f)] Local Poincare inequalities\\[-.3cm]
$$
2(t -s)\Gamma_t(P_{t,s} u)\le
P_{t,s}(u^2) - (P_{t,s}u)^2\le 2(t -s)P_{t,s}(\Gamma_s u)$$\\[-.8cm]
 \item[\rm g)] Local logarithmic Sobolev inequalities\\[-.3cm]
$$
(t -s)\frac{\Gamma_t(P_{t,s} u)}{P_{t,s} u}\le
P_{t,s}(u\,\log u) -  (P_{t,s}u)\log (P_{t,s}u)\le (t -s)P_{t,s}\Big(\frac{\Gamma_s u}{u}\Big)$$\\[-.8cm]
 \item[\rm h)] Dimension-free Harnack inequality: $\forall \alpha>1$\\[-.3cm]
$$(P_{t,s}u)^\alpha(y)\le P_{t,s}u^\alpha(x)\cdot \exp\Big(\frac{
\alpha \d_t^2(x,y)}{4(\alpha-1)(t-s)}\Big)$$\\[-.8cm]
\item[\rm i)] Log Harnack inequality: \\[-.3cm]
$$ P_{t,s}(\log u)(x)\le
\log P_{t,s}u(y)
+\frac{
 \d_t^2(x,y)}{4(t-s)}.$$
\end{enumerate}
\end{theorem}

With these concepts and results, a robust theory of super-Ricci 
flows is established --- being regarded as a parabolic analogue to singular spaces with lower Ricci bounds.
In the smooth case,  deeper insights and more powerful estimates 
require to restrict oneself to \emph{Ricci flows}
rather than  super-Ricci 
flows, see e.g.~\cite{haslhofer2015weak,kleiner2017singular,bamler2017uniqueness,kuwada2021monotonicity}.
To deal with similar questions in the singular case, 
first of all needs a  
synthetic notion of upper Ricci bounds, see next subsection.

For related current research on  lower Ricci bounds in time-like directions on Lorentzian manifolds and on
Einstein equation in general relativity, see \cite{mccann2018displacement,mondino2018optimal,cavalletti2020optimal}.

\subsection{Second Order Calculus, Upper Ricci Bounds, and Transformation Formulas}

So far, on $\RCD$-space we only dealt with the canonical $1^{st}$ order calculus for (real valued) functions on these spaces. The setting, however, allows us to go far beyond this.

\begin{theorem}[\cite{gigli2009second,gigli2018second,braun2020heat,gigli2020korevaar,gigli2020differential}] Given an $\RCD(K,\infty)$-space $(\X,\d,\m)$, there exist well established concepts of 
\begin{itemize}
\item[$\ast$] a powerful $2^{nd}$ order calculus on $\X$ including a consistent notion of  Ricci tensor (the lower bound of which coincides with the synthetic lower Ricci bound in terms of semiconvextiy of the entropy),  
\item[$\ast$] the heat flow on 1-forms on $\X$ which among others leads to the celebrated Hess-Schrader-Uhlenbrock inequality
$$|P_t \,\d f|\le e^{-Kt}\,P_t\,|\d f|,$$
\item[$\ast$] harmonic maps from $\X$ into metric spaces $(\sf Y,\d_{\sf Y})$, typically of nonpositive curvature, 
%coming with a Lipschitz regularity result.
based on Sobolev calculus and approximation of  energy densities for maps between metric spaces.
\end{itemize}
\end{theorem}

In a different direction, a challenging goal is to provide synthetic characterizations of  upper Ricci bounds $\Ric \le L$. Indeed, various of the (equivalent) synthetic characterizations of lower Ricci bounds admit partial converses. 
However, these converse characterizations are not necessarily equivalent to each other. Moreover, any such characterizations will certainly be not as powerful as the corresponding lower bound. 
Typically, the upper Ricci bounds are asymptotic estimates whereas the lower Ricci bounds are uniform estimates.

\begin{theorem}[\cite{SturmRemarks}]  Weak synthetic characterizations of  upper Ricci bounds for an $\RCD(K,N)$-space $(\X,\d,\m)$
\begin{itemize}
\item[$\ast$] in terms of partial $L$-concavity of the Boltzmann entropy and
\item[$\ast$] in terms of the heat kernel asymptotics
\end{itemize}
are equivalent to each other.

More precisely, a weak upper bound $L$ for the Ricci curvature is given by
 $$L:=\sup_{z}\limsup_{x,y\to z}\, \eta(x,y)$$  where for all $x,y\in X$,
\begin{eqnarray*}
\lefteqn{\eta(x,y):=}\\
&&=\lim_{\epsilon\to0}\  \inf\Big\{
\frac1{W_2^2(\rho^0,\rho^1)}\cdot\Big[ \partial^-_a S(\rho^a)\big|_{a=1}-\partial_a^+ S(\rho^a)\big|_{a=0}
\Big]: \ \big(\rho^a\big)_{a\in[0,1]} \mbox{ geodesic},\\
&&\qquad\qquad\qquad S(\rho^0)<\infty, \ S(\rho^1)<\infty, \  \supp[\rho^0]\subset B_\epsilon(x), \  \supp[\rho^1]\subset B_\epsilon(y)\Big\}\\
&&=\lim_{\epsilon\to0} \ \inf\Big\{
-\partial_t^+  \log W_2\big( P_t\mu,P_t\nu\big)\big|_{t=0}: \
\supp[\mu]\subset B_\epsilon(x), \supp[\nu]\subset B_\epsilon(y)\Big\}.
\end{eqnarray*}
\end{theorem}

\begin{remark} For weighted Riemannian manifolds $(\M,\g,e^{-f}d\mathrm{vol}_\g)$,
\begin{equation*}\label{sharp}
\ric_f(x,y)\le   \eta(x,y)
\le \ric_f(x,y)+\sigma(x,y)\cdot\tan^2\Big( \sqrt{\sigma(x,y)} \,\d(x,y)/2\Big)
\end{equation*}
provided $x$ and $y$ are not conjugate. Here
$\ric_f(x,y)=
 \int_{0}^1 \Ric_f(\dot\gamma^a, \dot\gamma^a)/ |\dot\gamma^a|^{2}\,da$
 denotes 
 the average Bakry--\'Emery--Ricci curvature along  the (unique) geodesic  $\gamma=(\gamma^a)_{a\in[0,1]}$   from $x$ to $y$,
and 
$\sigma(x,y)$ denotes the maximal modulus of the Riemannian curvature along this geodesic.
\end{remark}
Similar as other approaches (e.g.~\cite{naber2013characterizations}), these weak upper Ricci bounds will not be able to detect the positive Ricci curvature sitting in the tip of a cone over a circle of length $<2\pi$. A slightly stronger notion will detect it.
\begin{theorem}[\cite{ErbarSturm20}] If a metric cone has both sided (`strong') Ricci bounds $K$ and $L$ in the sense of $\RCD(K,\infty)$ and 
$$-\liminf_{x,y\to z}\liminf_{t\to0} \frac1{t}\log \frac{W_2\big(P_t\delta_x,P_t\delta_y\big)}{\d(x,y)}\le L\qquad(\forall z\in\X)$$
then it is the flat Euclidean space (of some integer dimension).
 \end{theorem}

A crucial property of the class of $\RCD$-spaces is that it is preserved under transformations of measure and metric of the underlying spaces, and that there exist explicit formulas for the transformation of the parameters $K$ and $N$ in the
curvature-dimension condition $\CD(K,N)$.

To be more specific, let an mm-space $(\X,\d,\m)$ be given as well as continuous (``weight'') functions $V,W$ on $\X$. In terms of them, define the transformed mm-space 
 $(\X,\d',\m')$ with $\m':=e^V \m$ and
\begin{equation*}
\d'(x,y)
:=\inf\bigg\{ \int_0^1 |\dot\gamma_t|\cdot e^{W(\gamma_t)}\,dt: \ \gamma:[0,1]\to X \mbox{ rectifiable, } \gamma_0=x, \gamma_1=y\bigg\}.
\end{equation*}
If $
\int |\nabla u|^2\,d\m$ on $L^2(\X,\m)
$
denotes the Dirichlet form (``Cheeger energy'') associated with $(\X,\d,\m)$, then the  Dirichlet form associated with the transformed mm-space is given by
$$
\int |\nabla u|^2 e^{V-2W}\,d\m\quad\mbox{on } L^2\big(\X,e^{V}\m\big).
$$

\begin{theorem}[\cite{sturm2018ricci,han2021curvature}] \label{time}
If $(\X,\d,\m)$ satisfies $\RCD(K,N)$ for finite $K,N\in\R$ and if $V,W\in W^{2,\infty}(X)$ then for each $N'>N$ there exists an explicitly given $K'$ s.t. $(\X,\d',\m')$ satisfies $\RCD(K',N')$. 

(If $W=0$ then also $N=N'=\infty$ is admissible; if $V=NW$ then also $N'=N$ is admissible.)
\end{theorem}
Let us illustrate this result in three special cases of particular importance:
\begin{itemize}
\item[$\ast$] $W=0$ (``Drift Transformation''):
$$K'=K-\sup_{f,x}\frac1{|\nabla f|^2}\Big[ \mathrm{Hess} V(\nabla f,\nabla f)+\frac1{N'-N}\langle\nabla V,\nabla f\rangle^2\Big](x)$$
\item[$\ast$] $V=2W$ (``Time Change''):
\begin{align*}K'
&=\inf_x e^{-2W}\big[K- \Delta W - \frac{[(N-2)(N'-2)]_+}{N'-N}\,|\nabla W|^2\big](x)
\end{align*}
\item[$\ast$] $V=NW$ (``Conformal Transformation''):\qquad $N'=N$ and
\begin{align*}K'=\inf_x e^{-2W}\bigg[K&-\Big[\Delta W+(N-2)|\nabla W|^2 \Big]\\
&-\sup_{f}\frac{N-2}{|\nabla f|^2}\Big[\mathrm{Hess} W(\nabla f,\nabla f)-\langle \nabla W,\nabla f\rangle^2 \Big]\bigg](x).
\end{align*}
\end{itemize}
The first of these cases is well-studied in the setting of Bakry--\'Emery calculus (and also in the setting of synthetic Ricci bounds for mm-spaces). It is the only case where also $N=\infty$ is admitted. The last of these cases is well-known in Riemannian geometry but has not been considered before in singular settings.
A particular feature of the second case is that the transformation formula for the Ricci bound only depends on bounds for $|\nabla W|$ and $\Delta W$ (and thus extends to distribution-valued Ricci bounds in case of $W\in\Lip(X)$, see next subsection).

\subsection{Distribution-valued Ricci Bounds}

Uniform lower Ricci bounds of the form $\CD(K,\infty)$ on mm-spaces
\begin{itemize}
\item
are preserved for Neumann Laplacian on convex subsets, but
\item
never hold
for Neumann Laplacian on non-convex subsets.
\end{itemize}
The goal thus is 
\begin{itemize}
\item
to find appropriate modification for  non-convex subsets
\item
to replace constant $K$, by function $k$, measure $\kappa$, distribution, \ldots
\end{itemize}

\begin{theorem}[\cite{BHS21}] \label{BHS} Given an infinitesimally Hilbertian mm-space $(\X,\d,\m)$ and  a
lower bounded, lower semicontinuous function $k:\X\to \R$, the following are equivalent

\begin{enumerate}
\item Curvature-dimension condition $\CD(k,\infty)$ with variable $k$: \ 
$\forall \mu_0,\mu_1 \in \Pz(X):\ \exists$\
$W_2$-geodesic $(\mu_t)_t=({e_t}_*\boldsymbol\nu)_t$\  s.t. $\forall t\in[0,1]$ with $g_{s,t}:=(1-s)t\wedge s(1-t)$:
$$\Ent(\mu_t)
 \leq  (1-t)\,\Ent(\mu_0) + t\,\Ent(\mu_1) -
 \int
 \int_0^1 k(\gamma_s)\, g_{s,t}\,ds\,|\dot\gamma|^2\,\boldsymbol\nu(d\gamma)
$$

\item Gradient estimate: 
$$\left|\nabla P_t u\right|(x)\le \mathbb E_x\Big[e^{-\int_0^t k(B_s)ds}\cdot |\nabla u|(B_t)\Big]$$

\item Bochner inequality 
$\BE_2(k,\infty)$: 
$$
\frac12\Delta|\nabla u|^2-\langle \nabla u, \nabla\Delta u\rangle
\ge k\cdot |\nabla u|^2$$

\item $\forall \mu_1,\mu_2:\ \exists$\ coupled pair of Brownian motions $(B^1_{t/2})_{t\ge0}$, $(B^2_{t/2})_{t\ge0}$ with given initial distributions such that  a.s. for all $s<t$\
 $$
 \d(B^1_t,B^2_t)\le e^{-\int_s^t \overline k(B^1_r,B^2_r)dr}\cdot \d(B^1_s,B^2_s)$$
with $\overline k(x_0,x_1):=\sup\big\{\int_0^1 k(\gamma_u)du: \ \gamma_0=x_0,\gamma_1=x_1, \ \gamma \text{ geodesic}\big\}$.

\end{enumerate}
\end{theorem}
For extensions to 
$(k,N)$-versions, see \cite{EKS,ketterer2017geometry,Sturm20distribution}.

To proceed towards distribution-valued Ricci bounds,
 define the spaces $W^{1,p}(\X)$ for $p\in[1,\infty]$, 
 put $W_*^{1,\infty}(\X):=\{f\in W^{1,2}_{loc}(\X): \| |\nabla f|\|_{L^\infty}<\infty\}$, and denote by $W^{-1,\infty}(\X)$  the topological dual of
  $$W^{1,1+}(\X):=\big\{f\in L^1(\X): \ f_n:=f\wedge n\vee(-n)\in W^{1,2}(\X), \ \sup_n \big\| |\nabla f_n|\big\|_{L^1}<\infty\big\}.$$

\begin{definition} Given  $\kappa\in W^{-1,\infty}(\X)$,
we say that the Bochner inequality $\BE_1(\kappa,\infty)$ holds iff $|\nabla f| \in W^{1,2}$ for all $f\in {\rm D}(\Delta)$, and
\begin{equation*}\label{H3-inequ}
-\int_X \langle \nabla|\nabla f|,\nabla\phi\rangle + \frac1{|\nabla f|}\langle\nabla f,\nabla \Delta f\rangle\phi\,dm \geq
\big\langle |\nabla f|\,\phi,\kappa\big\rangle_{W^{1,1}, W^{-1,\infty}}
\end{equation*}
for all $f\in {\rm D}(\Delta)$ with $\Delta f\in W^{1,2}$ and all  nonnegative $\phi \in W^{1,2}$.
\end{definition}

 Given $\kappa\in W^{-1,\infty}(X)$, 
we
define a closed, lower bounded  bilinear form ${\mathcal E}^\kappa$ on $L^2(X)$ by
$${\mathcal E}^\kappa(f,g):={\mathcal E}(f,g)+
\langle f\,g,\kappa\rangle_{W^{1,1+}, W^{-1,\infty}}
$$
for $f,g\in {\rm Dom}(\mathcal E^\kappa):=W^{1,2}(\X)$. 
Associated to it, there is a  strongly continuous, positivity preserving semigroup $(P^\kappa_t)_{t\ge0}$ on $L^2(\X)$.

\begin{theorem}[\cite{Sturm20distribution}] The Bochner inequality $\BE_1(\kappa,\infty)$ is equivalent to the 
gradient estimate
\begin{equation}
 |\nabla P_{t}f|
 \le P^{\kappa}_t\big( |\nabla f|\big).
\end{equation}
\end{theorem}

To gain a better understanding of the semigroup $(P^\kappa_t)_{t\ge0}$, assume that 
 $\kappa=-\underline\Delta\psi$ for some $\psi\in W^{1,\infty}$. 

\begin{theorem}[\cite{chen2002girsanov,Sturm20distribution}] 
Then 
\begin{equation}\mathcal E^\kappa(f,g)=\mathcal E(f,g)+\mathcal E(fg,\psi)\end{equation}
and
\begin{equation}P^\kappa_{t/2} f(x)={\mathbb E}_x\big[e^{N^\psi_t}f(B_t)\big]\end{equation}
where 
 $({\mathbb P}_x, (B_{t})_{t\ge0})$ denotes Brownian motion starting in $x\in X$,
and
$N^\psi$ is the zero energy part in the Fukushima decomposition, i.e.
$N^\psi_t=\psi(B_t)-\psi(B_0)-M^\psi_t$.
\end{theorem}

If $\psi\in {\rm Dom}(\Delta)$ then
$N^\psi_{t}=\frac12\int_0^t \Delta\psi(B_{s})ds$ --- in consistency with the previous Theorem \ref{BHS}.

\begin{remark} 
The concept of \emph{tamed spaces} proposed by Erbar--Rigoni--Sturm--Tamanini \cite{ERST} generalizes the previous approach to  distribution-valued lower Ricci bounds in various respects:
\begin{itemize}
\item[$\ast$] the objects under consideration are strongly local, quasi-regular Dirichlet spaces $(\X,\mathcal E,\m)$ (rather than infinitesimally Hilbertian mm-spaces $(\X,\d,\m)$);
\item[$\ast$] the Ricci bounds are formulated in terms of distributions $\kappa\in W^{-1,2}_{qloc}(\X)$ (rather than $\kappa\in W^{-1,\infty}(\X)$); for such distributions $\kappa$ which lie quasi locally in the dual of $W^{1,2}(\X)$, the previous ansatz for defining the semigroup $(P^\kappa_t)_{t>0}$ still works with appropriate sequences of localizing stopping times;
\item[$\ast$] in addition, the distributions $\kappa$ are assumed to be moderate in the sense that
$$\sup_{t\le 1, x\in\X}  P^\kappa_t1(x)<\infty.$$
This reminds of the Kato condition but is significantly more general since it does not require any decomposition of $\kappa$ into positive and negative parts. It  always holds if  $\kappa=-\underline\Delta\psi$ for some  $\psi\in {\rm Lip}_b(X)$.

\end{itemize}
\end{remark}

\begin{example} The prime examples of \emph{tamed spaces} are provided by:
\begin{enumerate}
\item ground state transformation of Hamiltonian for molecules \cite{guneysu2020molecules,braun2021heat}; yields curvature bounds in terms of unbounded functions in the Kato class;
\item Riemannian Lipschitz manifolds with lower Ricci bound in the Kato class \cite{rose2019li,carron2019geometric,carron2021limits};
\item time change of $\RCD(K,N)$-spaces with $W\in {\rm Lip}_b(\X)$ (cf.~Theorem \ref{time}); typically yields curvature bounds $\kappa$ which are not signed measures;

\item restriction of  $\RCD(K,N)$-spaces to (convex or non-convex) subsets $\Y\subset\X$ or, in other words, Laplacian with Neumann boundary conditions; yields curvature bounds in terms of signed measures  $\kappa=k\,m + \ell \,\sigma$, see below.
\end{enumerate}
\end{example}

Assume that $(\X,\d,\m)$ satisfies an $\RCD(k,N)$-condition with variable $k:\X\to \R$ and finite $N$.
Let a closed subset $\Y\subset \X$ be given which can be represented as sub-level set $\Y=\{V\le 0\}$ for some semiconvex function $V: \X\to\R$ with $|\nabla V|=1$ on 
$\partial\Y$.
Typically, $V$ is the signed distance functions
$V=\d(\,.\,, \Y)-\d(\,.\,,\X\setminus \Y)$.

 A function $\ell:\X\to \R$ is regarded as  ``generalized lower  bound for the curvature (or second fundamental) form of $\partial \Y$'' iff it is  a synthetic lower bound for the Hessian of $V$.
 
\begin{example} Assume that $\X$ is an Alexandrov space with  sectional 
curvature $\ge 0$ and that $\Y\subset \X$ satisfies an exterior ball condition: 
 $\forall z\in\partial \Y: \exists$ 
  ball $B_r(x)\subset \complement \Y$ with $z\in\partial B_r(x)$.
 Then 
 $\ell(z):=
-\frac1{r(z)}
$ is  a lower bound for the curvature  of  $\partial \Y$.
\end{example}

 Under weak  regularity assumptions, the distributional Laplacian $\sigma_\Y:=\underline\Delta V^+$ is a (nonnegative) measure which then will be regarded as ``the surface measure of $\partial \Y$''.

\begin{theorem}[\cite{Sturm20distribution}] Under weak regularity assumptions on $V$ and $\ell$, the restricted space $(\Y,\d_\Y,\m_\Y)$ satisfies a Bakry--\'Emery condition $\BE_1(\kappa,\infty)$ with a 
signed measure valued Ricci bound
\begin{equation}\kappa=k\cdot m_\Y +\ell\cdot\sigma_\Y.\end{equation}
Thus the Neumann heat semigroup on $\Y$ satisfies 
\begin{eqnarray}
| \nabla P_{t}^\Y u|(x)&\le
&\mathbb E_{x}\Big[ |\nabla u|(B_{t})\cdot e^{-\int_0^t  k(B_{s})ds}\cdot e^{-\int_0^t  \ell(B_{s})dL_s}\Big]
\end{eqnarray}
where  $(B_{s/2})_{s\ge0}$ denotes Brownian motion in $\Y$ and 
$(L_s)_{s\ge0}$ the continuous additive functional associated with $\sigma_\Y$.
\end{theorem}

For smooth subsets in Riemannian manifolds, this kind of gradient estimate 
--- with $(L_s)_{s\ge0}$ being the \emph{local time} of the boundary ---
has been firstly derived by Hsu \cite{hsu2002multiplicative}, cf.~also \cite{wang2009second,cheng2017functional}.

Let us illustrate the power of the above estimates with two simple examples:  the ball and its complement.

\begin{corollary}
Let $(\X,\d,\m)$ be an $N$-dimensional Alexandrov space ($N\ge 3$) with  $\Ric \ge-1$ and ${\sf sec}\le0$.
Then for  $\Y:=\X\setminus B_r(z)$,
\begin{eqnarray*}
\big| \nabla P_{t/2}^\Y f\big|(x)&\le
&\mathbb E^\Y_{x}\Big[ e^{t/2+\frac1{2r}L^{\partial \Y}_t}
\cdot \big|\nabla f(B^\Y_{t})\big|\Big].
\end{eqnarray*}
In particular, $\mathrm{Lip}( P_{t/2}^\Y f)\le \sup_x\mathbb E^\Y_{x}\big[ e^{t/2+\frac1{2r}L^{\partial \Y}_t}
\big]\cdot  \mathrm{Lip}(f)$ and 
\begin{eqnarray}
\big| \nabla P_{t/2}^Yf\big|^2(x)\le 
 e^{Ct+C'\sqrt{t}}\cdot {P_{t/2}^\Y\big| \nabla f\big|^2(x)}.
\end{eqnarray}
\end{corollary}

Upper and lower bound of curvature (here 0 and -1, resp) can be chosen to be any numbers.
Note that \emph{no estimate} of the form 
\begin{eqnarray*}
{\big| \nabla P_{t/2}^\Y f\big|^2(x)}\le e^{Ct}\cdot {P_{t/2}^\Y\big| \nabla f\big|^2(x)}
\end{eqnarray*}
 can hold true due to the non-convexity of $\Y$. Thus it is \emph{necessary} to take into account the singular contribution arising from the negative curvature of the boundary. 
 
 In the next example, the singular contribution arising from the positive curvature of the boundary can be ignored. However, taking it into account will significantly \emph{improve} the gradient estimate.

\begin{corollary}
Let $(\X,\d,\m)$ be an $N$-dimensional Alexandrov space with $\Ric \ge0$ and ${\sf sec} \le 1$. Then for $\Y:=\overline B_r(z)$ for some $z\in \X$ and $r\in (0,\pi/4)$.
\begin{eqnarray*}
\big| \nabla P_{t/2}^Yf\big|(x)&\le
&\mathbb E^\Y_{x}\Big[ e^{-\frac{\cot r}{2} \, L^{\partial \Y}_t}
\cdot \big|\nabla f(B^Y_{t})\big|\Big].
\end{eqnarray*}
In particular, $\mathrm{Lip}( P_{t/2}^Yf)\le \sup_x\mathbb E^\Y_{x}\big[ e^{-\frac{\cot r}{2} \, L^{\partial \Y}_t}\cdot  \mathrm{Lip}(f)
\big]$ and 
\begin{eqnarray}
{\big| \nabla P_{t/2}^\Y f\big|^2(x)}
\le e^{-t\frac{N-1}{2}\cot^2r+1}\cdot {P_{t/2}^\Y\big| \nabla f\big|^2(x)}.
\end{eqnarray}
\end{corollary}

Taking into account the curvature of the boundary,
allows us to derive  a positive lower bound for the spectral gap (without involving any diameter bound and despite possibly vanishing Ricci curvature  in the interior).
\begin{corollary} In the previous setting,
\begin{equation}\lambda_1\ge \frac{N-1}{2}\cot^2r.\end{equation}
\end{corollary}

\subsection{Synthetic Ricci Bounds \quad --- \quad Extended Settings}

In order to summarize recent developments concerning  synthetic Ricci bounds for singular spaces, let us recall the previously presented 
\begin{itemize}
\item[a)] \emph{heat flow on time-dependent mm-spaces and super-Ricci flows,}
\item[b)] 
\emph{second order calculus, upper Ricci bounds, and transformation formulas,}
\item[c)] \emph{distribution-valued lower Ricci bounds,} 
\end{itemize}
and then move on to further developments in extended settings
\begin{enumerate}

\item[d)] \emph{discrete mm-spaces:}

For discrete mm-spaces $(\X,\d,\m)$, the synthetic Ricci bounds as introduced above will be meaningless since there will be no non-constant geodesics w.r.t.~the Kantorovich--Wasserstein metric $W_2$ as defined in \eqref{W2}. 
This disadvantage can be overcome by resorting to a modified Kantorovich--Wasserstein metric based on a subtle discrete version of the Benamou--Brenier formula. This way, the heat flow can again be characterized as the gradient flow of the entropy ~\cite{Maas11,Mielke11}. 

And synthetic Ricci bounds defined in terms of semiconvexity of the entropy w.r.t.~this modified metric are intimately linked to equilibration properties of the heat flow, see e.g.~\cite{erbar2012ricci,gozlan2014displacement,erbar2015discrete,erbar2018poincare,erbar2019entropic}.
Challenging questions address homogenization \cite{gigli2013gromov,gladbach2020scaling,gladbach2021homogenisation} and evolution under curvature flows \cite{erbar2020super}.
Related --- but in general different --- concepts of synthetic Ricci bounds are based on discrete versions of the Bakry--\'Emery condition, see e.g.~\cite{bauer2011ollivier,liu2018bakry,weber2017characterizing,fathi2018curvature,cushing2020rigidity}.
\item[e)] 
\emph{non-commutative spaces:}

Inspired by the synthetic Ricci bounds for discrete spaces, an analogous concept also has been proposed for non-commutative spaces, with remarkable insights e.g.~for (ergodic) quantum Markov semigroups on tracial or finite-dimensional unital $C^*$-algebras, in particular, equilibration rate estimates for the fermionic Ornstein--Uhlenbeck semigroup and for
Bose Ornstein-Uhlenbeck semigroups
\cite{mielke2013dissipative,carlen2014analog,carlen2020non,hornshaw2019quantum,wirth2021dual}.

\item[f)] 
\emph{Dirichlet boundary conditions:}

For a long time it seemed that OT techniques could not be used to analyze the heat flow with Dirichlet boundary conditions. Only recently, Profeta--Sturm \cite{ProfetaSturm} overcame the problem of mass absorption by considering \emph{charged particles} (which are either particles or anti-particles),
and this way succeeded in finding a characterization for the heat flow as a gradient flow for the entropy.
Passing from particles to charged particles technically corresponds to passing from a space $\X$ to its \emph{doubling}.
Functional inequalities for the Dirichlet heat flow thus are closely linked to those for the doubled space.
For recent progress concerning the challenging  \emph{problem of gluing
convex subsets in $\RCD$-spaces}, see \cite{KaKeSt}.

\end{enumerate}

\bibliographystyle{amsalpha}
\bibliography{biblio-ecm-sturm}

\end{document}